\documentclass{article}
\usepackage{amssymb}
\usepackage{amsmath}

\input{tcilatex}
\begin{document}

\author{Piotr WILCZEK}
\title{Some Representation Theorem for nonreflexive Banach space \\
ultrapowers under the Continuum Hypothesis.}
\date{}
\maketitle

\begin{abstract}
In this paper it will be shown that assuming the \textit{Continuum Hypothesis%
} $(CH)$ every nonreflexive Banach space ultrapower $E^{I}/\mathcal{U}$ is
isometrically isomorphic to the space of continuous, bounded and real-valued
functions on the Stone-Cech remainder $\omega ^{\ast }$. This Representation
Theorem will be helpful in proving some facts from geometry and topology of
nonreflexive Banach space ultrapowers.
\end{abstract}


$2010$\textit{\ Mathematical Subject Classification}. $46B08,46B20$, $46B25$%
.\bigskip 

$\mathbf{1.}$\textbf{\ Introduction\bigskip }

Although the ultrapower construction was initially developed within \textit{%
Model Theory} it exerted the great impact on almost all other branches of
mathematics (especially such as \textit{Algebra} and \textit{Set Theory}).
In the field rendered as \textit{Banach Space Theory} this construction was
introduced in the mid $60$s by Bretagnolle, Dacunha-Castelle and Krivine $($%
\cite{BDCK, DCK}$)$. Recall that a Banach space $E$ is said to be finite
dimensional if and only if its unit ball is compact, i.e., if and only if
for every bounded family $(x_{i})_{i\in I}$ and for every ultrafilter $%
\mathcal{U}$ on the set $I$ the so-called $\mathcal{U-}$limit%
\begin{equation*}
\underset{i,\mathcal{U}}{\lim }x_{i}
\end{equation*}%
exists. But if a Banach space $E$ is infinite dimensional, then it is
possible to enlarge $E$ to a Banach space $\widehat{E}$ by adjoining to
every bounded family $(x_{i})_{i\in I}$ in $E$ an element $\widehat{x}\in 
\widehat{E}$ such that $\left\Vert \widehat{x}\right\Vert =\underset{i,%
\mathcal{U}}{\lim }\left\Vert x_{i}\right\Vert $. This construction is
termed \textit{Banach space ultrapower}. Suppose that $\left( E_{i}\right)
_{i\in I}$ is an index family of Banach spaces. Then define%
\begin{equation*}
\ell _{\infty }\left( E_{i}\right) =\left\{ \left( x_{i}\right) :x_{i}\in
E_{i}\text{ and }\left\Vert \left( x_{i}\right) \right\Vert _{\infty
}<\infty \right\} \text{.}
\end{equation*}%
It is easily observed that $\ell _{\infty }\left( E_{i}\right) $ is the
Banach space of all bounded families $(x_{i})\in \underset{i\in I}{\dprod }%
E_{i}$ endowed with the norm given by $\left\Vert (x_{i})\right\Vert
_{\infty }=\underset{i\in I}{\sup }\left\Vert x_{i}\right\Vert _{E_{i}}$. If 
$\mathcal{U}$ is an ultrafilter on the index set $I$, then it is always
possible to determine $\underset{i,\mathcal{U}}{\lim }\left\Vert
x_{i}\right\Vert _{E_{i}}$. Then it is seen that $\mathcal{N((}x_{i}\mathcal{%
))=}\underset{i,\mathcal{U}}{\lim }\left\Vert x_{i}\right\Vert _{E_{i}}$ is
a seminorm on $\ell _{\infty }(E_{i})$. Consequently, the kernel of $%
\mathcal{N}$ is given by 
\begin{equation*}
\mathcal{N}_{\mathcal{U}}=\left\{ x=(x_{i})\in \ell _{\infty }(E_{i}):%
\underset{i,\mathcal{U}}{\lim }\left\Vert x_{i}\right\Vert =0\right\} \text{.%
}
\end{equation*}%
It follows that $\mathcal{N}_{\mathcal{U}}$ is a closed ideal in the Banach
space $\ell _{\infty }(E_{i})$. Thus it is possible to define the quotient
space of the form:%
\begin{equation*}
\ell _{\infty }\left( E_{i}\right) /\mathcal{N}_{\mathcal{U}}\text{.}
\end{equation*}%
This quotient is said to be the \textit{ultraproduct} of the family of
Banach spaces $\left( E_{i}\right) _{i\in I}$. If $E_{i}=E$ for every $i\in
I $, then the space $\ell _{\infty }(E_{i})/\mathcal{U}$ is called the 
\textit{ultrapower} of $E$ and is symbolized by $\ell _{\infty }(E)/\mathcal{%
N}_{\mathcal{U}}$ $($or by $E^{I}/\mathcal{U})$. Therefore, it can be
observed that - from the model-theoretical point of view - the above
construction can be considered as the consequence of eliminating the
elements of infinite norm from an ordinary $($i.e., algebraic$)$ ultrapower
and dividing it by infinitesimal $($\cite{AK, BBHU, BDCK, DCK, H, HI}$)$. \
\ \ \ \ \ \ \ \ \ \ \ \ \ \ \ \ \ \ \ \ \ \ \ \ \ \ \ \ \ \ \ \ \ \ \ \ \ \
\ \ \ \ \ \ \ \ \ \ \ \ \ \ \ \ \ \ \ \ \ \ \ \ \ \ \ \ \ \ \ \ \ \ \ \ \ \
\ \ \ \ \ \ \ \ \ \ 

Recall that the famous \textit{\L o\'{s} Theorem} - also known as the 
\textit{Fundamental Theorem on Ultraproducts} - asserting that any
first-order formula is true in the ultraproduct $\mathcal{M=}\underset{i\in I%
}{\dprod }M_{i}/\mathcal{U}$ $($where $\mathcal{M}$ is any first order
structure$)$ if and only if the set of indices $i$ such that the formula is
true in $M_{i}$ is a member of $\mathcal{U}$ can be easily adapted to the
case of the so-called \textit{positive bounded formulas} which are more
adequate for considering metric structures $($cf. Proposition $9.2$ in \cite%
{HI}$)$. \ \ \ \ \ \ \ \ \ \ \ \ \ \ \ \ \ \ \ \ \ \ \ \ \ \ \ \ \ \ \ \ \ \
\ \ \ \ \ \ \ \ \ \ \ \ \ \ \ \ \ \ \ \ \ \ \ \ \ \ \ \ \ \ \ \ \ \ \ \ \ \
\ \ \ \ \ \ \ \ \ \ \ \ \ \ \ \ \ \ \ \ \ \ \ \ \ \ \ \ \ \ \ \ \ \ \ \ \ \
\ \ \ \ \ \ \ \ \ \ \ \ \ \ \ \ \ \ \ \ \ \ 

If $(x_{i})$ is a family in $\ell _{\infty }(E)/\mathcal{N}_{\mathcal{U}}$,
then let us denote by $(x_{i})_{\mathcal{U}}$ the \textit{equivalence class}
of $(x_{i})$ in the ultrapower $\ell _{\infty }(E)/\mathcal{N}_{\mathcal{U}}$%
. If $E$ is any Banach space and $\ell _{\infty }(E)/\mathcal{N}_{\mathcal{U}%
}$ is its ultrapower, then the mapping $x\rightarrow (x_{i})_{\mathcal{U}}$,
where $x_{i}=x$ for every $i\in I$, constitutes an \textit{isometry} of $E$
into $\ell _{\infty }(E)/\mathcal{N}_{\mathcal{U}}$. Consequently, it
becomes obvious that $E$ is a subspace of $\ell _{\infty }(E)/\mathcal{N}_{%
\mathcal{U}}$. The above mentioned isometric embedding generally is not 
\textit{onto}. On the other hand, the above isometry is surjective if the
ultrafilter $\mathcal{U}$ is principal or the space $E$ is finite
dimensional $($\cite{HI}$)$. \ \ \ \ \ \ \ \ \ \ \ \ \ \ \ \ \ \ \ \ \ \ \ \
\ \ \ \ \ \ \ \ \ \ \ \ \ \ \ \ \ \ \ \ \ \ \ \ \ \ \ \ \ \ \ \ \ \ \ \ \ \
\ \ \ \ \ \ \ \ \ \ \ \ \ \ \ \ \ \ \ \ \ \ \ \ \ \ \ \ \ \ \ \ \ \ \ \ \ \
\ \ \ \ \ \ \ \ \ \ \ \ \ \ \ \ \ \ \ \ \ \ \ \ \ \ \ \ \ \ \ \ \ \ \ \ \ \
\ \ \ \ \ \ \ \ \ \ \ \ \ \ \ \ \ \ \ \ \ \ \ \ \ \ \ \ \ \ \ \ \ \ \ \ \ \
\ \ \ \ \ \ \ \ \ \ \ \ \ \ \ \ \ \ \ \ \ \ \ \ \ \ \ \ \ \ \ \ \ \ \ \ \ \
\ \ \ \ \ \ \ \ \ \ \ \ \ \ \ \ \ \ \ \ \ \ \ \ \ \ \ \ \ \ \ \ \ \ \ \ \ \
\ \ \ \ \ \ \ \ \ \ \ \ \ \ \ \ \ \ \ \ \ \ \ \ \ \ \ \ \ \ \ \ \ \ \ \ \ \
\ \ \ \ \ \ \ \ \ \ \ 

Recall that an index family $\left( x_{i}\right) _{i\in I}$ in any
topological space converges to the point $x$ with respect to an ultrafilter $%
\mathcal{U}$, i.e., $\underset{i,\mathcal{U}}{\lim }=x$ if for every open
set $V$ containing the point $x$ it follows that the set $\left\{ i\in
I:x_{i}\in V\right\} $ belongs to $\mathcal{U}$. \ \ \ \ \ \ \ \ \ \ \ \ \ \
\ \ \ \ \ \ \ \ \ \ \ \ \ \ \ \ \ \ \ \ \ \ \ \ \ \ \ \ \ \ \ \ \ \ \ \ \ \
\ \ \ \ \ \ \ \ \ \ \ \ \ \ \ \ \ \ \ \ \ \ \ \ \ \ \ \ \ \ \ \ \ \ \ \ \ \
\ \ \ \ \ \ \ \ \ \ \ \ \ \ \ \ \ \ \ \ \ \ \ \ \ \ \ \ \ \ \ \ \ \ \ \ \ \
\ \ \ \ \ \ \ \ \ \ \ \ \ \ \ \ \ \ \ \ \ \ \ \ \ \ \ \ \ \ \ \ \ \ \ \ \ \
\ \ \ \ \ \ \ \ \ \ \ \ \ \ \ \ \ \ \ \ \ \ \ \ \ \ \ \ \ \ \ \ \ \ \ \ \ \
\ \ \ \ \ \ \ \ \ \ \ \ \ \ \ \ \ \ \ \ \ \ \ \ \ \ \ \ \ \ \ \ \ \ \ \ \ \
\ \ \ \ \ \ \ \ \ \ \ \ 

It is known that structural questions about Banach space ultrapowers are
really interesting only when it is assumed that the considered ultrafilter $%
\mathcal{U}$ is \textit{countably incomplete}, i.e., it is possible to
single out a sequence $\left( U_{n}\right) $ of members of $\mathcal{U}$
such that $\underset{n}{\dbigcap }U_{n}=\varnothing $. This requirement
holds especially for free ultrafilters on $\mathbb{\omega }$ $($\cite{AK,
AK1, H}$)$.\ \ \ \ \ \ \ \ \ \ \ \ \ \ \ \ \ \ \ \ \ \ \ \ \ \ \ \ \ \ \ \ \
\ \ \ \ \ \ \ \ \ \ \ \ \ \ \ \ \ \ \ \ \ \ \ \ \ \ \ \ \ \ \ \ \ \ \ \ \ \
\ \ 

A Banach space $E$ is called superreflexive if and only if each ultrapower $%
\ell _{\infty }(E)/\mathcal{N}_{\mathcal{U}}$ is reflexive $($\cite{AK, H}$)$%
.\ \ \ \ \ \ \ \ \ \ \ \ \ \ \ \ \ \ \ \ \ \ \ \ \ \ \ \ \ \ \ \ \ \ \ \ \ \
\ \ \ \ \ \ \ \ \ \ \ \ \ \ \ \ \ \ \ \ \ \ \ \ \ \ \ \ \ \ \ \ \ \ \ \ \ \
\ \ \ \ \ \ \ \ \ \ \ \ \ \ \ \ \ \ \ \ \ \ \ \ \ \ \ \ \ \ \ \ \ \ \ \ \ \
\ \ \ \ \ \ \ \ \ \ \ \ \ \ \ \ \ \ \ \ \ \ \ \ \ \ \ \ \ \ \ \ \ \ \ \ \ \
\ \ \ \ \ \ \ \ \ \ \ \ \ \ \ \ \ \ \ \ \ \ \ \ \ \ \ \ \ \ \ \ \ \ \ \ \ \
\ \ \ \ \ \ \ \ \ \ \ \ \ \ \ \ \ \ \ \ \ \ \ \ \ \ \ \ \ \ \ \ \ \ \ \ \ \
\ \ \ \ \ \ \ \ \ \ \ \ \ \ \ \ \ \ \ \ \ \ \ \ \ \ \ \ \ \ \ \ \ \ \ \ \ \
\ \ \ \ \ \ \ \ \ \ \ \ \ \ \ \ \ \ \ \ \ \ \ \ \ \ \ \ \ \ \ \ \ \ \ 

In this study we denote by $\omega $ the set of all natural numbers, i.e., $%
\omega =\left\{ 1,2,...\right\} $ and its \textit{Stone-Cech compactification%
} by $\beta \omega $. Then the so-called \textit{Stone-Cech remainder} of $%
\beta \omega $, i.e., the space $\beta \omega \backslash \omega $ is denoted
by $\omega ^{\ast }$ $($\cite{W}$)$. The space $\beta \omega $ can be
identified with the set of ultrafilters on $\omega $ under the topology
generated by the sets of the form $\left\{ F:U\in F\right\} $ for all sets $%
U\subseteq \omega $. On the other hand, it should be mentioned that the set $%
\omega $ corresponds to the set of principal ultrafilters and the Stone-Cech
remainder $\omega ^{\ast }$ corresponds to the set of free ultrafilters on $%
\omega $ $($\cite{W}$)$.\ \ \ \ \ \ \ \ \ \ \ \ \ \ \ \ \ \ \ \ \ \ \ \ \ \
\ \ \ \ \ \ \ \ \ \ \ \ \ \ \ \ \ \ \ \ \ \ \ \ \ \ \ \ \ \ \ \ \ \ \ \ \ \
\ \ \ \ \ \ \ \ \ \ \ \ \ \ \ \ \ \ \ \ \ \ \ \ \ \ \ \ \ \ \ \ \ \ \ \ \ \
\ \ \ \ \ \ \ \ \ \ \ \ \ \ \ \ \ \ \ \ \ \ \ \ \ \ \ \ \ \ \ \ \ \ \ \ \ \
\ \ \ \ \ \ \ \ \ \ \ \ \ \ \ \ \ \ \ \ \ \ \ \ \ \ \ \ \ \ \ \ \ \ \ \ \ \
\ \ \ \ \ \ \ \ \ \ \ \ \ \ \ \ \ \ \ \ \ \ \ \ \ \ \ \ \ \ \ \ \ \ \ \ \ \
\ \ \ \ \ \ \ \ \ \ \ \ \ \ \ \ \ \ \ \ \ \ \ \ \ \ \ \ \ \ \ \ \ \ \ \ \ \
\ \ \ \ \ \ \ \ \ \ \ \ \ \ \ \ \ \ \ \ \ \ \ \ \ \ \ \ \ \ \ \ \ \ \ \ \ \
\ \ \ \ \ \ \ \ \ \ \ \ \ \ \ \ From \textit{Set-Theoretical Topology} it is
known that the so-called \textit{Parovicenko space} $X$ is identified with a
topological space satisfying the following conditions: $1)$ $X$ is compact
and Hausdorff, $2)$ $X$ has no isolated points, $3)$ $X$ has the weight $%
\mathfrak{c}$, $4)$ every nonempty $G_{\delta }$ subset of $X$ has nonempty
interior, $5)$ every two disjoint open $F_{\sigma }$ subsets of $X$ have
disjoint closures. In $1963$ I. I. Parovi\v{c}enko proved that assuming the
continuum hypothesis $(CH)$ every Parovicenko space $X$ is isomorphic to $%
\omega ^{\ast }$ $($\cite{W}$)$. \ \ \ \ \ \ \ \ \ \ \ \ \ \ \ \ \ \ \ \ \ \
\ \ \ \ \ \ \ \ \ \ \ \ \ \ \ \ \ \ \ \ \ \ \ \ \ \ \ \ \ \ \ \ \ \ \ \ \ \
\ \ \ \ \ \ \ \ \ \ \ \ \ \ \ \ \ \ \ \ \ \ \ \ \ \ \ \ \ \ \ \ \ \ \ \ \ \
\ \ \ \ \ \ \ \ \ \ \ \ \ \ \ \ \ \ \ \ \ \ \ \ \ \ \ \ \ \ \ \ \ \ \ \ \ \
\ \ \ \ \ \ \ \ \ \ \ \ \ \ \ \ \ \ \ \ \ \ \ \ \ \ \ \ \ \ \ \ \ \ \ \ \ \
\ \ \ \ \ \ \ \ \ \ \ \ \ \ \ \ \ \ \ \ \ \ \ \ \ \ \ \ \ \ \ \ \ \ \ \ \ \
\ \ \ \ \ \ 

Since Banach space ultrapowers were introduced into the field of \textit{%
Functional Analysis} a considerable numbers of paper employing this
methodology can be observed $($\cite{AK, BBHU, BDCK, DCK, H, HI} and papers
cited there$)$. But on the other hand, it should be noted that there exists
relatively little papers concerning the topological and geometrical
structure of these model-theoretical objects. Consequently, it is hoped that
the following article will be useful in solving the problems concerning
these spaces.\ \ \ \ \ \ \ \ \ \ \ \ \ \ \ \ \ \ \ \ \ \ \ \ \ \ \ \ \ \ \ \
\ \ \ \ \ \ \ \ \ \ \ \ \ \ \ \ \ \ \ \ \ \ \ \ \ \ \ \ \ \ \ \ \ \ \ \ \ \
\ \ \ \ \ \ \ \ \ \ \ \ \ \ \ \ \ \ \ \ \ \ \ \ \ \ \ \ \ \ \ \ \ \ \ \ \ \
\ \ \ \ \ \ \ \ \ \ \ \ \ \ \ \ \ \ \ \ \ \ \ \ \ \ \ \ \ \ \ \ \ \ \ \ \ \
\ \ \ \ \ \ \ \ \ \ \ \ \ \ \ \ \ \ \ \ \ \ \ \ \ \ \ \ \ \ \ \ \ \ \ \ \ \
\ \ \ \ \ \ \ \ \ \ \ \ \ \ \ \ \ \ \ \ \ \ \ \ \ \ \ \ \ \ \ \ \ \ \ \ \ \
\ \ \ \ \ \ \ \ \ \ \ \ \ \ \ \ \ \ \ \ \ \ \ \ \ \ \ \ \ \ \ \ \ \ \ \ \ \
\ \ \ \ \ \ \ \ \ \ \ \ \ \ \ \ \ \ \ \ \ \ \ \ \ \ \ \ \ \ \ \ \ \ \ \ \ \
\ \ \ \ \ \ \ \ \ \ \ \ \ \ \ \ \ \ \ \ \ \ \ 

In our studies it will be shown that under $CH$ every Banach space
ultrapower can be alternatively represented in the form of the space of
continuous, bounded and real-valued functions defined on the Parovicenko
space $\omega ^{\ast }$. Namely, assuming $CH$ the \textit{Representation
Theorem} for nonreflexive Banach space ultrapowers will be obtained. This
new result ascertains that if $CH$ holds and $E$ is any infinite dimensional
nonsuperreflexive Banach space, then the Banach space ultrapower $\ell
_{\infty }(E)/\mathcal{N}_{\mathcal{U}}$ is isometrically isomorphic to the
space of continuous, bounded and real-valued functions on the Stone-Cech
remainder $\omega ^{\ast }$. The congruence $\ell _{\infty }(E)/\mathcal{N}_{%
\mathcal{U}}\cong C(\omega ^{\ast })$ which holds under $CH$ is of central
importance in our further studies. This isometric isomorphism will enable to
prove that nonreflexive Banach space ultrapowers are never dual spaces, the
unit ball of any nonreflexive Banach space ultrapower has an abundance of
extreme points and no smooth points. Also the structure of complemented
subspaces of $\ell _{\infty }(E)/\mathcal{N}_{\mathcal{U}}$ will be
elucidated.

In this place it should be mentioned that our results are obtained under the
assumption that all considered infinite dimensional Banach spaces are
nonsuperreflexive and - consequently - their ultrapowers are nonreflexive.
It is unknown if the condition of nonsuperreflexivity can be weakened $($or
modified$)$ in order to formulate our Representation Theorem.\bigskip 

$\mathbf{2.}$\textbf{\ Notation\bigskip }

\textit{Dual} of any Banach space $E$ is denoted by $E^{\ast }$. The \textit{%
unit ball} of $E$ is symbolized by $B_{E}$ and the \textit{unit sphere} of $%
E $ by $S_{E}$. Denote by $A$ any convex subset in a Banach space $E$. Then
the \textit{convex hull} of $A$, denoted by $co(A)$, is identified with the
smallest convex set containing $A$. The \textit{closed convex hull} of $A$,
denoted by $\overline{co}(A)$, is the smallest closed convex set which
contains the subset $A$. If $A$ is convex, then any point $x\in A$ is said
to be an \textit{extreme point} of $A$ if whenever $x=\lambda
x_{1}+(1-\lambda )x_{2}$ for $0<\lambda <1$, then $x=x_{1}=x_{2}$.
Alternatively, the point $x$ is an extreme point of the subset $A$ if $%
A\backslash \{x\}$ is still convex. Denote by $\partial _{e}(A)$ the set of
all extreme points of the subset $A$. On the other hand, the \textit{smooth
points} of the unit ball of the space $C(T)$, where $T$ is compact and
Hausdorff, are identified with the functions $f\in C(T)$ such that $%
\left\Vert f\right\Vert =\sup \left\{ \left\vert f(t)\right\vert :t\in
T\right\} =1$. This means that these functions peak at some $t_{0}\in T$, i.
e., $\left\vert f(t_{0})\right\vert =1>\left\vert f(t)\right\vert $ for all $%
t\in T$ such that $t\neq t_{0}$. If some function $f$ peaks at isolated
points of the space $T$, then $f$ is said to be the\textit{\ point of Fr\'{e}%
chet differentiability} of the supremum norm on the space $C(T)$ $($\cite%
{AK1}$)$. \ \ \ \ \ \ \ \ \ \ \ \ \ \ \ \ \ \ \ \ \ \ \ \ \ \ \ \ \ \ \ \ \
\ \ \ \ \ \ \ \ \ \ \ \ \ \ \ \ \ \ \ \ \ \ \ \ \ \ \ \ \ \ \ \ \ \ \ \ \ \
\ \ \ \ \ \ \ \ \ \ \ \ \ \ \ \ \ \ \ \ \ \ \ \ \ \ \ \ \ \ \ \ \ \ \ \ \ \
\ \ \ \ \ \ \ \ \ \ \ \ \ \ \ \ \ \ \ \ \ \ \ \ \ \ \ \ \ \ \ \ \ \ \ \ \ \
\ \ \ \ \ \ \ \ \ \ \ \ \ \ \ \ \ \ \ \ \ \ \ \ \ \ \ \ \ \ \ \ \ \ \ \ \ \
\ \ \ \ \ \ \ \ \ \ \ \ \ \ \ \ \ \ \ \ \ \ \ \ \ \ \ \ \ \ \ \ \ \ \ \ \ \
\ \ \ \ \ \ \ \ \ \ \ \ \ \ \ \ \ \ \ \ \ \ \ \ \ The symbol $\cong $ is
used in order to denote the relation of isometric isomorphism between Banach
spaces. \ \ \ \ \ \ \ \ \ \ \ \ \ \ \ \ \ \ \ \ \ \ \ \ \ \ \ \ \ \ \ \ \ \
\ \ \ \ \ \ \ \ \ \ \ \ \ \ \ \ \ \ \ \ \ \ \ \ \ \ \ \ \ \ \ \ \ \ \ \ \ \
\ \ \ \ \ \ \ \ \ \ \ \ \ \ \ \ \ \ \ \ \ \ \ \ \ \ \ \ \ \ \ \ \ \ \ \ \ \
\ \ \ \ \ \ \ \ \ \ \ \ \ \ \ \ \ \ \ \ \ \ \ \ \ \ \ \ \ \ \ \ \ \ \ \ \ \
\ \ \ \ 

If $E$ and $F$ are two Banach spaces, then an operator $Q:X\rightarrow Y$ is
said to be \textit{compact} $($\textit{weakly compact}, respectively$)$ if
the closure of $Q(B_{E})$ is compact $($weakly compact, respectively$)$ $($%
\cite{AK1}$)$.\ \ \ \ \bigskip\ \ \ \ \ \ \ \ \ \ \ \ \ \ \ \ \ \ \ \ \ \ \
\ \ \ \ \ \ \ \ \ \ \ \ \ \ \ \ \ \ \ \ \ \ \ \ \ \ \ \ \ \ \ \ \ \ \ \ \ \
\ \ \ \ \ \ \ \ \ \ \ \ \ \ \ \ \ \ \ \ \ \ \ \ \ \ \ \ \ \ \ \ \ \ \ \ \ \
\ \ \ \ \ \ \ \ \ \ \ \ \ \ \ \ \ \ \ \ \ \ \ \ \ \ \ \ \ \ \ \ \ \ \ \ \ \
\ \ \ \ \ \ \ \ \ \ \ \ \ \ \ \ \ \ \ \ \ \ \ \ \ \ \ \ \ \ \ \ \ \ \ \ \ \
\ \ \ \ \ \ \ \ \ \ \ \ \ \ \ \ \ \ \ \ \ \ \ \ \ \ \ \ \ \ \ \ \ \ \ \ \ 

$\mathbf{3.}$\textbf{\ The Representation Theorem for nonreflexive Banach
space ultrapowers \bigskip }

Identifying any Tychonoff space $T$ with completely regular and Hausdorff
space it can be shown that its Stone-Cech compactification $\beta T$ can be
represented in the following form. Suppose that $C(T)$ is the Banach space
of all continuous, bounded and real-valued functions on $T$ with the norm
defined by $\left\Vert f\right\Vert =\sup \left\{ \left\vert f(t)\right\vert
:t\in T\right\} $ and assume that $B_{C(T)^{\ast }}$ denotes the closed unit
ball of the dual space $C(T)^{\ast }$. If we identify every element $t\in T$
with the evaluation functional $\phi _{t}\in B_{C(T)^{\ast }}$, where $\phi
_{t}(f)=f(t)$ for $f\in C(T)$, then it is possible to represent $\beta T$ as
the weak$^{\ast }$-closure of the set $\left\{ \phi _{t}:t\in T\right\} $ in 
$B_{C(T)^{\ast }}$ and $T$ can be understood as a dense subset of $\beta T$.
Consequently, every function $f\in C(T)$ has a unique norm-preserving
extension $\widehat{f}\in C(\beta T)$ $($cf. \cite{LW, M}$)$. \ \ \ \ \ \ \
\ \ \ \ \ \ \ \ \ \ \ \ \ \ \ \ \ \ \ \ \ \ \ \ \ \ \ \ \ \ \ \ \ \ \ \ \ \
\ \ \ \ \ \ \ \ \ \ \ \ \ \ \ \ \ \ \ \ \ \ \ \ \ \ \ \ \ \ \ \ \ \ \ \ \ \
\ \ \ \ \ \ \ \ \ \ \ \ \ \ \ \ \ \ \ \ \ \ \ \ \ \ \ \ \ \ \ \ \ \ \ \ \ \
\ \ \ \ \ \ \ \ \ \ \ \ \ \ \ \ \ \ \ \ \ \ \ \ \ \ \ \ \ \ \ \ \ \ \ \ \ \
\ \ \ \ \ \ \ \ \ \ \ \ \ \ \ \ \ \ \ \ \ \ \ \ \ \ \ \ \ \ \ \ \ \ \ \ \ \
\ \ \ \ \ \ \ \ \ \ \ \ \ \ \ \ \ \ \ \ \ \ \ \ \ \ \ \ \ \ \ \ \ \ \ \ \ \
\ \ \ \ \ \ \ \ \ \ \ \ \ \ \ \ \ \ \ \ \ \ \ \ \ \ \ \ If the set of all
natural numbers $\omega $ has the discrete topology and $E$ is any Banach
space, then $\ell _{\infty }(E)=C(\omega )$. If we define the restriction
mapping $R:C(\beta \omega )\rightarrow C(\omega )$ by $R(\widehat{f})=%
\widehat{f}\upharpoonright \omega $ for each $\widehat{f}\in C(\beta \omega )
$, then it can be concluded that $R$ is a linear isometry of $C(\beta \omega
)$ onto $C(\omega )$. Hence, the following proposition can be
asserted:\medskip 

\textbf{Proposition }$\mathbf{1}$ $(CH)$. \textit{Let }$E$\textit{\ be any
infinite dimensional nonsuperreflexive Banach space and let }$\ell _{\infty
}\left( E\right) $\textit{\ be the }$\ell _{\infty }-$\textit{sum of
countably many copies of }$E$\textit{. Then }$\ell _{\infty }\left( E\right) 
$\textit{\ is isometrically isomorphic to the space }$C(\beta \omega )$%
\textit{. Symbolically}%
\begin{equation*}
\ell _{\infty }(E)\cong C(\beta \omega )\text{.}
\end{equation*}%
\medskip Now suppose that $I$ is the closed ideal in the space $C(\beta
\omega )$ consisting of functions which vanish on $\omega ^{\ast }$, i.e., $%
I=\left\{ \widehat{f}\in C(\beta \omega ):\widehat{f}(t)=0\text{ for all }%
t\in \omega ^{\ast }\right\} $. Next, assume that the space $c_{0}(\omega )$
consists of functions in $C(\omega )$ which vanish at infinity, i.e., $%
c_{0}(\omega )=\left\{ f\in C(\omega ):\text{ for each }\varepsilon >0\text{%
, }\left\{ t\in \omega :\left\vert f(t)\right\vert >\varepsilon \right\} 
\text{ is finite}\right\} $. On the other hand, it is known that if $%
\mathcal{U}$ is a nontrivial ultrafilter on $\mathbb{\omega }$ and the
sequence $(x_{n})$ converges to the point $x$ in the topology of the space $E
$, then $\left( x_{n}\right) $ converges to $x$ with respect to $\mathcal{U}$%
, i.e., $\underset{\mathcal{U}}{\lim }x_{n}=x$. This follows from the simple
observation that if $V$ denotes any neighborhood of $x$, then the set $%
\left\{ i:x_{i}\notin V\right\} $ is finite and the nontriviality of $%
\mathcal{U}$ implies that the set $\left\{ i:x_{i}\in V\right\} $ belongs to 
$\mathcal{U}$ $($cf. Proposition $2.2$ in \cite{AK}$)$. Basing on this facts
it can be claimed that $c_{0}(\omega )=\mathcal{N}_{\mathcal{U}}$. Also the
restriction mapping $R:I\rightarrow \mathcal{N}_{\mathcal{U}}$ defines a
linear isometry from $I$ onto $\mathcal{N}_{\mathcal{U}}$. Consequently, we
arrive at the following proposition:\medskip 

\textbf{Proposition} $\mathbf{2}$ $(CH)$. \textit{Let }$E$\textit{\ be any
infinite dimensional \textit{nonsuperreflexive}e Banach space, }$N_{\mathcal{%
U}}=\left\{ \left( x_{i}\right) \in \ell _{\infty }(E):\underset{\mathcal{U}}%
{\lim }\left\Vert x_{i}\right\Vert =0\right\} $\textit{\ be the closed ideal
in the }$\ell _{\infty }-$\textit{sum of countably many copies of }$E$%
\textit{\ and }

$I=\left\{ \widehat{f}\in C(\beta \omega ):\widehat{f}(t)=0\text{ for all }%
t\in \omega ^{\ast }\right\} $\textit{\ be the closed ideal in the space }$%
C(\beta \omega )$\textit{. Then both ideals are isometrically isomorphic.
Symbolically}%
\begin{equation*}
\mathcal{N}_{\mathcal{U}}\cong I\text{.}
\end{equation*}%
\medskip Further, assume that the mapping $\sigma :C(\beta \omega
)/I\rightarrow C(\omega ^{\ast })$ defined by $\sigma \left( \widehat{f}%
+I\right) =\widehat{f}\upharpoonright \omega ^{\ast }$ for each function $%
\widehat{f}\in C(\beta \omega )$ constitutes a linear isometry from $C(\beta
\omega )/I$ onto $C(\omega ^{\ast })$. Then the following corollary can be
easily obtained:\medskip 

\textbf{Corollary }$\mathbf{3}$ $(CH)$.\textit{\ Let }$E$\textit{\ be any
infinite dimensional \textit{nonsuperreflexive}e Banach space and }$\ell
_{\infty }(E)/N_{\mathcal{U}}$\textit{\ be its ultrapower. Then the Banach
space ultrapower is isometrically isomorphic to the space of continuous,
bounded and real-valued functions defined on the Stone-Cech remainder }$%
C(\omega ^{\ast })$\textit{, i.e., }%
\begin{equation*}
\ell _{\infty }(E)/\mathcal{N}_{\mathcal{U}}\cong C(\omega ^{\ast })\text{.}
\end{equation*}%
\medskip 

Corollary $3$ can be regarded as the \textit{Representation Theorem} since
it asserts that - under $CH$- each Banach space ultrapower $\ell _{\infty
}(E)/\mathcal{N}_{\mathcal{U}}$ $($where $E$ is any infinite dimensional
nonsuperreflexive\textit{\ }Banach space$)$ can be isometrically isomorphic
represented in the form of the space of continuous, bounded and real-valued
functions on the Parovicenko space $C(\omega ^{\ast })$.\bigskip 

$\mathbf{4.}$\textbf{\ Nonreflexive Banach space ultrapowers are never dual
spaces\bigskip }

Suppose that $T$ is a compact Hausdorff space. Then $\mathcal{B}$ denotes
the $\sigma -$algebra of Borel subsets of $T$ and $rca(T,\mathcal{B})$ is
the Banach space of regular and countably additive Borel measures $\mu $ on $%
T$ endowed with bounded variation. This norm is given by the variation $\mu $
on $T$, i.e., $\left\Vert \mu \right\Vert =\left\vert \mu \right\vert
(T)=\sup \overset{n}{\underset{i=1}{\dsum }}\left\vert \mu \left(
A_{i}\right) \right\vert $ where the supremum ranges over all finite
partitions $\left\{ A_{1},A_{2},...,A_{n}\right\} $ of the space $T$. Now,
define in $rca(T,\mathcal{B})$ the norm closed proper cone containing
positive normal measures. Denote this cone by $N^{+}\left( T,\mathcal{B}%
\right) $. The measure $\mu $ is said to be normal if $\mu (B)=0$ for each
Borel set $B$ which is meager in $T$. Also assume that $N^{+}\left( T,%
\mathcal{B}\right) $ generates the closed ideal in the space $rca\left( T,%
\mathcal{B}\right) $ which is denoted by $N\left( T,\mathcal{B}\right) $.
Every measure $\mu $ in the space $rca\left( T,\mathcal{B}\right) $ is
supported on the set of the form $S(\mu )=\dbigcap \left\{ F\subseteq T:F%
\text{ is closed and }\left\vert \mu \right\vert (F)=\left\vert \mu
\right\vert (T)\right\} $. Recall that the compact Hausdorff space is said
to be \textit{hyperstonian} if $T$ is extremally disconnected and the sum $%
\dbigcup \left\{ S(\mu ):\mu \in N^{+}\left( T,\mathcal{B}\right) \right\} $
constitutes the dense set in $T$. A. Grothendieck in $($\cite{Gr}$)$ \
proved that any compact Hausdorff space $T$ is congruent to a dual space if $%
T$ is hyperstonian. Namely, the following theorem can be formulated:\medskip 

\textbf{Theorem }$\mathbf{4}$ $($Grothendieck$)$. \textit{If }$T$\textit{\
is a compact Hausdorff space and }$X$\textit{\ is any Banach space, then }$%
L:C(T)\rightarrow X^{\ast }$\textit{\ is an isometric isomorphism of }$C(T)$%
\textit{\ onto }$X^{\ast }$\textit{\ and }$J:X\rightarrow X^{\ast \ast }$%
\textit{\ is the canonical embedding, then }

$i)$\textit{\ }$T$\textit{\ is hyperstonian,}

$ii)$\textit{\ }$L^{\ast }\circ J$\textit{\ is an isometric isomorphism of }$%
X$\textit{\ onto }$N(T,B)$.\medskip 

Also in $($\cite{D}$)$ the converse of this theorem was proved. Namely, it
was demonstrated that if $T$ is hyperstonian, then $N(T,\mathcal{B})^{\ast }$
is congruent to $C(T)$. These considerations can be easily adapted to the
case of Banach space ultrapowers. The following conclusion can be
obtained:\medskip 

\textbf{Theorem }$\mathbf{5}$. \textit{Let }$E$\textit{\ be any infinite
dimensional nonsuperreflexive Banach space and }$\ell _{\infty }(E)/\mathcal{%
N}_{\mathcal{U}}$ \textit{be its ultrapower. Then its ultrapower }$\ell
_{\infty }(E)/N_{\mathcal{U}}$\textit{\ is not a dual space.}

\textit{Proof}. From our Representation Theorem for nonreflexive Banach
space ultrapowers it follows that $\ell _{\infty }(E)/\mathcal{N}_{\mathcal{U%
}}\cong C(\omega ^{\ast })$. From this identification and from the fact that 
$\omega ^{\ast }$ is not extremally disconnected $($\cite{W}$)$ it is
straightforward to see that the Banach space $\ell _{\infty }(E)/\mathcal{N}%
_{\mathcal{U}}$ is not a dual space. $\square \bigskip $

$\mathbf{5.}$\textbf{\ Geometry of nonreflexive Banach space ultrapower }$%
\ell _{\infty }(E)/\mathcal{N}_{\mathcal{U}}\bigskip $

$\mathbf{5.1}$ \textbf{Extreme points of} \textbf{the unit ball} $B_{\ell
_{\infty }(E)/\mathcal{N}_{\mathcal{U}}}$.\medskip 

It will be shown that the unit ball $B_{\ell _{\infty }(E)/\mathcal{N}_{%
\mathcal{U}}}$ $($where $E$ is any infinite dimensional nonsuperreflexive
Banach space$)$ has an abundance of extreme points. Even it can be proved
that $B_{\ell _{\infty }(E)/\mathcal{N}_{\mathcal{U}}}=\overline{co}\left(
\partial _{e}\left( B_{\ell _{\infty }(E)/\mathcal{N}_{\mathcal{U}}}\right)
\right) $. \ \ \ \ \ \ \ \ \ \ \ \ \ \ \ \ \ \ \ \ \ \ \ \ \ \ \ \ \ \ \ \ \
\ \ \ \ \ \ \ \ \ \ \ \ \ \ \ \ \ \ \ \ \ \ \ \ \ \ \ \ \ \ \ \ \ \ \ \ \ \
\ \ \ \ \ \ \ \ \ \ \ \ \ \ \ \ \ \ \ \ \ \ \ \ \ \ \ \ \ \ \ \ \ \ \ \ \ \
\ \ \ \ \ \ \ \ \ \ \ \ \ \ \ \ \ \ \ \ \ \ \ \ \ \ \ \ \ \ \ \ \ \ \ \ \ \
\ \ \ \ \ \ \ \ \ \ \ \ \ \ \ \ \ \ \ \ \ \ \ \ \ \ \ \ \ \ \ \ \ \ \ \ \ \
\ \ \ \ \ \ \ \ \ \ \ \ \ \ \ \ \ \ \ \ \ \ \ \ \ \ \ \ \ \ \ \ \ \ \ \ \ \
\ \ \ \ \ \ \ \ \ \ \ \ \ \ \ \ \ \ \ \ \ \ \ \ \ \ \ \ \ \ \ \ \ \ \ \ \ \
\ \ \ Firstly, it will be demonstrated that every extreme point of the unit
ball $B_{\ell _{\infty }(E)/\mathcal{N}_{\mathcal{U}}}$ can be represented
as the image $($with respect to some quotient mapping$)$ of an extreme point
in $B_{\ell _{\infty }(E)}$ $($cf. \cite{LW}$)$.\medskip 

\textbf{Theorem }$\mathbf{6}$ $\mathbf{(}CH\mathbf{).}$ \textit{Let }$E$%
\textit{\ be any infinite dimensional nonsuperreflexive Banach space. If }$%
q:\ell _{\infty }(E)\rightarrow \ell _{\infty }(E)/N_{\mathcal{U}}$\textit{\
is the quotient mapping, then }$\partial _{e}(B_{\ell _{\infty }(E)/\mathcal{%
N}_{\mathcal{U}}})=q(\partial _{e}(B_{\ell _{\infty }(E)}))$\textit{.}

\textit{Proof.} From the Representation Theorem it is known that the spaces

$\ell _{\infty }(E)/\mathcal{N}_{\mathcal{U}}$ and $C(\omega ^{\ast })$ are
congruent. Consequently, it must be prove that $\partial _{e}[B_{C(\beta
\omega ^{\ast })}]=\pi \lbrack \partial _{e}(B_{C(\beta \omega )})]$ $($or
isometrically isomorphic: $\partial _{e}[B_{\ell _{\infty }(E)/\mathcal{N}_{%
\mathcal{U}}}]=\pi \lbrack \partial _{e}(B_{\ell _{\infty }(E)})])$ where $%
\pi :C(\beta \omega )\rightarrow C(\omega ^{\ast })$ is identified with the
quotient mapping defined by $\pi (f)=f\upharpoonright \omega ^{\ast }$ for
every $f\in C(\beta \omega )$. Recall that if $T$ is a compact Hausdorff
space, then $\left\vert f(t)\right\vert =1$ for all $t\in T$. Suppose that $%
\widehat{p}$ is an extreme point of $B_{C(\beta \omega ^{\ast })}$. Then $%
\left\vert \widehat{p}(t)\right\vert =1$ for all $t\in \omega ^{\ast }$. It
is possible for each $t\in \omega ^{\ast }$ to single out an open
neighborhood $V_{t}$ of $t$ in $\omega ^{\ast }$ such that $\widehat{p}%
\upharpoonright V_{t}$ is constant. Also it is possible to find out an open
neighborhood $\widehat{V}_{t}$ of $t$ in $\beta \omega $ such that $V_{t}=%
\widehat{V}_{t}\cap (\omega ^{\ast })$. From the fact that $\beta \omega
=\dbigcup \left\{ \widehat{V}_{t}:t\in \omega ^{\ast }\}\cup \{\{n\}:n\in
\omega \right\} $ is compact it can be deduced that there exists a finite
family $\left\{ \widehat{V}_{t_{1}},...,\widehat{V}_{t_{k}},\{n_{1}\},...,%
\{n_{j}\}\right\} $ covering $\beta \omega $. Also it can be claimed that $%
n_{i}\notin \widehat{V}_{t_{l}}$ for any $i,l$. Then it is possible to
introduce the mapping $p:\beta \omega \rightarrow \mathbb{R}$ defined by:%
\begin{eqnarray*}
p(t) &=&\widehat{p}\upharpoonright V_{t_{i}}\text{ for }t\in \widehat{V}%
_{t_{i}}\text{, }i=1,...,k\text{ } \\
&&\text{or} \\
p(t) &=&1\text{ for }t=n_{i}\text{, }i=1,...,j\text{.}
\end{eqnarray*}%
It can be observed that the mapping $p$ is continuous. Namely, if $t_{0}\in
\omega $, then $\left\{ t\right\} $ is an open neighborhood of $t_{0}$. If $%
t_{0}\in \omega ^{\ast }$, then $t_{0}\in V_{t_{i}}\subseteq \widehat{V}%
_{t_{i}}$ for some $i$ such that $1\leq i\leq k$. If it is assumed that the
net $t_{\delta }\rightarrow t_{0}$ is in $\beta \omega $, then $(t_{\delta })
$ is eventually in $\widehat{V}_{t_{i}}$ and eventually%
\begin{equation*}
p(t_{\delta })=\widehat{p}\upharpoonright V_{t_{i}}(t_{\delta })=\widehat{p}%
\upharpoonright V_{t_{i}}(t_{0})=p(t_{0})_{.}
\end{equation*}%
Therefore, it follows that $p(t_{\delta })\rightarrow p(t_{0})$ and $p$ is
continuous at $t_{0}.$ From the fact that $p:\beta \omega \rightarrow 
\mathbb{R}$ it is obvious that $p\in C(\beta \omega )$ and $\left\vert
p(t)\right\vert =1$ for all $t\in \beta \omega $ and - consequently - $p\in
\partial _{e}[B_{C(\beta \omega )}]$. It can be also observed that $\widehat{%
p}=p\upharpoonright \omega ^{\ast }=\pi (p)$. Consequently, it is
straightforward that $\partial _{e}[B_{C(\beta \omega \backslash \omega
)}]\subseteq \pi (\partial _{e}[B_{C(\beta \omega )}])$. Undoubtedly, this
inclusion can be reversed and the theorem is proved. $\square \medskip $

Recall that a closed subspace $M$ of a real normed linear space $E$ is said
to be \textit{proximinal} in $E$ if for each point $x\in E$ it is possible
to find out the point $y\in M$ such that $\left\Vert x-y\right\Vert =\inf
\left\{ \left\Vert x-z\right\Vert :z\in M\right\} $. Then it is possible to
state the following theorem $($cf. \cite{LW}$)$:\medskip 

\textbf{Theorem }$\mathbf{7}$ $(CH)$. \textit{Let }$E$\textit{\ be any
infinite dimensional nonsuperreflexive Banach space, }$\ell _{\infty }(E)$%
\textit{\ be the }$\ell _{\infty }-$\textit{sum of countably many copies of }%
$E$\textit{\ and }$N_{\mathcal{U}}=\left\{ \left( x_{i}\right) \in \ell
_{\infty }(E):\underset{\mathcal{U}}{\lim }\left\Vert x_{i}\right\Vert
=0\right\} $\textit{\ be the closed ideal in }$\ell _{\infty }(E)$\textit{.
Then the ideal }$N_{\mathcal{U}}$\textit{\ is proximinal in }$\ell _{\infty
}(E)$.

\textit{Proof.} Recall that $\mathcal{N}_{\mathcal{U}}$ is isometrically
isomorphic to the ideal $I=\{f\in C(\beta \omega ):f(t)=0$ for all $t\in
\omega ^{\ast }\}$ and $\ell _{\infty }(E)$ is congruent to the space $%
C(\beta \omega )$ $($i.e., $\mathcal{N}_{\mathcal{U}}\cong I$ and $\ell
_{\infty }(E)\cong C(\beta \omega )$, respectively$)$ $($Propositions $1$
and $2)$. Then it must be shown that the closed ideal $I$ is proximinal in
the space $C(\beta \omega )$. Suppose that $f\in C(\beta \omega )$ and $%
F=f\upharpoonright \omega ^{\ast }\in C(\omega ^{\ast })$. From Tietze's
Extension Theorem it is possible to find out the function $h\in C(\beta
\omega )$ such that $h\upharpoonright \omega ^{\ast }=F=f\upharpoonright
\omega ^{\ast }$ and $\left\Vert h\right\Vert =\left\Vert F\right\Vert =%
\underset{t\in \omega ^{\ast }}{\sup }\left\vert h(t)\right\vert $. It is
clear that $f-h\in I$ and $\left\Vert f+I\right\Vert =\left\Vert
h+I\right\Vert =\inf \{\left\Vert h-g\right\Vert :g\in I\}\leq \left\Vert
h\right\Vert =\left\Vert h\upharpoonright \omega ^{\ast }\right\Vert
=\left\Vert f\upharpoonright \omega ^{\ast }\right\Vert $. Also for any $%
g\in I$ it follows that $\left\Vert h-g\right\Vert \geq \left\Vert
(h-g)\upharpoonright \omega ^{\ast }\right\Vert =\left\Vert h\upharpoonright
\omega ^{\ast }\right\Vert =\left\Vert f\upharpoonright \omega ^{\ast
}\right\Vert $. Consequently, $\left\Vert f\upharpoonright \omega ^{\ast
}\right\Vert =\left\Vert f+I\right\Vert $. Suppose that $g_{0}=f-h$, then $%
g_{0}\in I$ and $\left\Vert f-g_{0}\right\Vert =\left\Vert h\right\Vert
=\left\Vert f\upharpoonright \omega ^{\ast }\right\Vert =\left\Vert
f+I\right\Vert $ and $\left\Vert f-g_{0}\right\Vert =\inf \{\left\Vert
f-g\right\Vert :g\in I\}$. Therefore, $I$ is proximinal in $C(\beta \omega )$
and from the identification of $\mathcal{N}_{\mathcal{U}}$ with $I$ it
follows that the ideal $\mathcal{N}_{\mathcal{U}}$ is proximinal in $\ell
_{\infty }(E)$. $\square \medskip $

In order to prove our main corollary we are forced to refer to the theorem
obtained by Godini $($cf. \cite{G}$)$. Namely:\medskip 

\textbf{Theorem }$\mathbf{8}$ $($Godini$)$. \textit{If }$E$\textit{\ is any
real normed linear space, }$M\subseteq E$\textit{\ is a closed subspace and }%
$r:E\rightarrow E/M$\textit{\ is the quotient mapping, then the following
conditions are equivalent:}

$1/$\textit{\ }$r(B_{E})=B_{E/M}$\textit{,}

$2/$\textit{\ }$r(B_{E})$\textit{\ is closed in }$E/M$\textit{,}

$3/$\textit{\ }$M$\textit{\ is proximinal in }$E$\textit{.\medskip }

Then it is possible to state the following corollary:\medskip 

\textbf{Corollary }$\mathbf{9}$\textbf{\ }$(CH)$. \textit{Let }$E$\textit{\
be any infinite dimensional nonsuperreflexive Banach space and }$\ell
_{\infty }(E)/N_{\mathcal{U}}$\textit{\ be its ultrapower. Then }

$B_{\ell _{\infty }(E)/\mathcal{N}_{\mathcal{U}}}=\overline{co}\left(
\partial _{e}\left( B_{\ell _{\infty }(E)/\mathcal{N}_{\mathcal{U}}}\right)
\right) $.

\textit{Proof}. It can be immediately seen that $B_{\ell _{\infty }(E)/%
\mathcal{N}_{\mathcal{U}}}=S_{B_{\ell _{\infty }(E)}}=$

$r(\overline{co}(\partial _{e}(B_{\ell _{\infty }(E)})))=\overline{co}%
(r(\partial _{e}(B_{\ell _{\infty }(E)})))=\overline{co}(\partial
_{e}(B_{\ell _{\infty }(E)/\mathcal{N}_{\mathcal{U}}}))$. $\square \medskip $

$\mathbf{5.2}$\textbf{\ Smooth points of} \textbf{the unit ball} $B_{\ell
_{\infty }(E)/\mathcal{N}_{\mathcal{U}}}$.\medskip 

In the previous section it was indicated that the quotient mapping $q:\ell
_{\infty }(E)\rightarrow \ell _{\infty }(E)/\mathcal{N}_{\mathcal{U}}$ for
any infinite dimensional\textit{\ }nonsuperreflexive Banach space $E$
preserves extreme points of the unit ball \textbf{\ }$B_{\ell _{\infty }(E)/%
\mathcal{N}_{\mathcal{U}}}$. But this does not hold for smooth points. It
can be demonstrated that the unit ball $B_{\ell _{\infty }(E)/\mathcal{N}_{%
\mathcal{U}}}$ has no smooth points $($cf. \cite{LW}$)$.\medskip 

\textbf{Theorem }$\mathbf{10}$ $(CH)$. \textit{Let }$E$\textit{\ be any
infinite dimensional nonsuperreflexive Banach space and }$\ell _{\infty
}(E)/N_{\mathcal{U}}$\textit{\ be its ultrapower. Then the unit ball }$%
B_{\ell _{\infty }(E)/\mathcal{N}_{\mathcal{U}}}$\textit{\ has no smooth
points.}

\textit{Proof. }From the identification $($under $CH)$ of Banach space
ultrapower $\ell _{\infty }(E)/\mathcal{N}_{\mathcal{U}}$ with the space $%
C(\omega ^{\ast })$ it follows that it must be demonstrated that for the
function $f\in C(\omega ^{\ast })$ such that $\left\Vert f\right\Vert =1$
the set $A=\{t\in \omega ^{\ast }:\left\vert f(t)\right\vert =\left\Vert
f\right\Vert =1\}$ is nonempty and contains more than one element. It can be
observed that the set $A=\overset{\infty }{\underset{n=1}{\dbigcap }}\left\{
t\in \omega ^{\ast }:\left\vert f(t)\right\vert >\left\Vert f\right\Vert -%
\frac{1}{n}\right\} $ is $G_{\delta }$ subset of $\omega ^{\ast }$. From $($%
\cite{W}$)$ it is known that $A$ has nonempty interior and - consequently -
contains a clopen subset of the Stone-Cech remainder $\omega ^{\ast }$. From
the fact that $\omega ^{\ast }$ has no isolated points it is deducible that $%
card(A)\geq 2$. Then it is obvious that the function $f$ is not a smooth
point of the unit ball $B_{C(\omega ^{\ast })}$ and our theorem is proved. $%
\square \medskip $

It should be observed that the smooth points of the unit ball $B_{\ell
_{\infty }(E)}\cong B_{C(\beta \omega )}$ can be identified with the points
of Fr\'{e}chet differentiability of the norm of $\ell _{\infty }(E)$. It
follows from the fact that if $f\in B_{C(\beta \omega )}$, $\left\Vert
f\right\Vert =1$ and $f$ peaks at $t_{0}\in \beta \omega $, then $t_{0}\in
\omega $ is an isolated point of the Stone-Cech compactification $\beta
\omega $.\bigskip 

$\mathbf{7.}$ \textbf{Complemented subspaces of }$\ell _{\infty }(E)/%
\mathcal{N}_{\mathcal{U}}\bigskip $

It is known that any closed subspace $M$ of any Banach space $E$ is said to
be \textit{complemented} in $E$ if the Banach space $E$ can be written as a
direct sum of $M$ and a closed subspace $N$ of $E$. Then the projection
operator $P:E\rightarrow M$ $($i.e., the mapping of $E$ onto $M$ along $N)$
is continuous; $M$ and $N$ are termed \textit{complementary subspaces} and $E
$ can be written as follows: $E=M\oplus N$. Rosenthal $($cf. \cite{DU}$)$
showed that if the space $T$ is extremally disconnected and a Banach space $E
$ has no copy of the space $\ell _{\infty }$, then every bounded linear
operator $U:C(T)\rightarrow E$ is weakly compact. It will be seen that
Rosenthal's result can be easily generalized to the case of nonreflexive
Banach space ultrapowers $\ell _{\infty }(E)/\mathcal{N}_{\mathcal{U}}$ $($%
cf. \cite{LW, M}$)$.\medskip 

\textbf{Theorem }$\mathbf{11}$ $(CH)$. \textit{Let }$E$\textit{\ be any
infinite dimensional nonsuperreflexive Banach space. If }$M$\textit{\ is an
infinite dimensional complemented subspace of a Banach space ultrapower }$%
\ell _{\infty }(E)/\mathcal{N}_{\mathcal{U}}$,\textit{\ then }$M$ \textit{%
contains a subspace which is isometrically isomorphic to the} $\ell _{\infty
}-$\textit{sum of countably many Banach spaces.}

\textit{Proof.} In this proof it must be demonstrated that if the subspace $M
$ contains no $\ell _{\infty }-$sum of countably many Banach spaces $E$,
then $M$ is finite dimensional. Namely, define the continuous projection
operator $P:\ell _{\infty }(E)/\mathcal{N}_{\mathcal{U}}\rightarrow M$ $($%
i.e., the mapping of $\ell _{\infty }(E)/\mathcal{N}_{\mathcal{U}}$ onto $M$
along its complement$)$ and the quotient mapping $q:\ell _{\infty
}(E)\rightarrow \ell _{\infty }(E)/\mathcal{N}_{\mathcal{U}}$. From the
above mentioned Rosenthal's result and the facts that $\ell _{\infty
}(E)\cong C(\beta \omega )$ and the space $\beta \omega $ is extremally
disconnected it follows that the composition operator $U=P\circ q$ is weakly
compact. Consequently, it can be seen that $U(B_{\ell _{\infty }(E)})$ is
relatively weakly compact and from Theorem $3$ it is deducible that $%
P(B_{\ell _{\infty }(E)/\mathcal{N}_{\mathcal{U}}})=U(B_{\ell _{\infty }(E)})
$ is also weakly compact. Hence, we obtain that the projection operator $P$
is weakly compact and - basing on the relation of congruence $\ell _{\infty
}(E)/\mathcal{N}_{\mathcal{U}}\cong C(\omega ^{\ast })$ - it is observed
that $P^{2}=P$ is compact. Then it is straightforward to see that the
subspace $M$ is finite dimensional. $\square \medskip $

Immediately we get the following result $($cf. \cite{LW, M}$)$.\medskip 

\textbf{Corollary} $\mathbf{12}$ $(CH)$. \textit{Let }$E$\textit{\ be any
infinite dimensional nonsuperreflexive Banach space and} $\ell _{\infty }(E)/%
\mathcal{N}_{\mathcal{U}}$ \textit{be} \textit{its ultrapower. Then} $\ell
_{\infty }(E)/\mathcal{N}_{\mathcal{U}}$ \textit{has no infinite dimensional
complemented subspaces which are separable or reflexive}.\medskip 

In the next theorem we are going to show that any nonreflexive Banach space
ultrapower $\ell _{\infty }(E)/\mathcal{N}_{\mathcal{U}}$ is isometrically
isomorphic to its square $\ell _{\infty }(E)/\mathcal{N}_{\mathcal{U}}\times
\ell _{\infty }(E)/\mathcal{N}_{\mathcal{U}}$ with an adequate norm $($cf. 
\cite{LW, M}$)$.\medskip 

\textbf{Theorem }$\mathbf{13}$ $(CH)$. \textit{Let }$E$\textit{\ be any
infinite dimensional nonsuperreflexive Banach space and} $\ell _{\infty }(E)/%
\mathcal{N}_{\mathcal{U}}$ \textit{be} \textit{its ultrapower. If the square
ultrapower }$\ell _{\infty }(E)/\mathcal{N}_{\mathcal{U}}\times \ell
_{\infty }(E)/\mathcal{N}_{\mathcal{U}}$ \textit{has the norm} $\left\Vert
(x,y)\right\Vert _{0}=\max (\left\Vert x\right\Vert ,\left\Vert y\right\Vert
)$ \textit{where} $x,y\in \ell _{\infty }(E)/\mathcal{N}_{\mathcal{U}}$, 
\textit{then there exists an isometric isomorphism} $Z$ of $\ell _{\infty
}(E)/\mathcal{N}_{\mathcal{U}}\times \ell _{\infty }(E)/\mathcal{N}_{%
\mathcal{U}}$ \textit{onto} $\ell _{\infty }(E)/\mathcal{N}_{\mathcal{U}}$.

\textit{Proof.} Let $A$ and $B$ be two nonempty disjoint clopen subsets of
the Stone-Cech remainder $\omega ^{\ast }$ such that $A\cup B=\omega ^{\ast }
$. Recall that the Parovicenco space $\omega ^{\ast }$ is totally
disconnected and each nonempty clopen subset of $\omega ^{\ast }$ has the
following form: $cl_{\beta \omega }D\backslash \omega $ where $D\subseteq
\omega $ is infinite. Also each such clopen subset is homeomorphic to $%
\omega ^{\ast }$. Then it is possible to find out the following
homeomorphisms: $\phi :A\rightarrow \omega ^{\ast }$ and $\phi :B\rightarrow
\omega ^{\ast }$. For two functions $f,g\in C(\omega ^{\ast })$ define the
operator $Q(f,g)=h$ such that%
\begin{eqnarray*}
h(t) &=&f(\varphi (t))\text{ if }t\in A\text{ } \\
&&\text{or} \\
h(t) &=&g(\varphi (t))\text{ if }t\in B\text{.}
\end{eqnarray*}%
It is obvious that $h\in C(\omega ^{\ast })$, for assume that $t_{\delta
}\rightarrow t$ and $t_{\delta },t\in \omega ^{\ast }$. If $t\in A$, then $%
t_{\delta }$ is eventually in $A$ and the sequence $t_{\delta }$ is
convergent to $t$ in $A$. As the result of this presupposition it is
obtained that $\varphi (t_{\delta })\rightarrow \varphi (t)$ and $%
h(t_{\delta })=(f\circ \varphi )(t_{\delta })\rightarrow (f\circ \varphi
)(t)=h(t)$. Analogously, if $t\in B$, then $h(t_{\delta })\rightarrow h(t)$.
Thus $h\in C(\omega ^{\ast })$ and the operator $Q:C(\omega ^{\ast })\times
C(\omega ^{\ast })\rightarrow C(\omega ^{\ast })$ is linear. Now we want to
show that $Q$ is onto. Namely, for $h\in C(\omega ^{\ast })$ suppose that $%
f=h\circ \varphi ^{-1}$ and $g=h\circ \varphi ^{-1}$. Then it is obtained
that $Q(f,g)=h$ and 
\begin{eqnarray*}
\left\Vert Q(f,g)\right\Vert  &=&\sup \left\{ \left\vert
Q(f,g)(t)\right\vert :t\in \beta \omega \backslash \omega \right\}  \\
&=&\max \left( \sup \left\{ \left\vert f(\varphi (t))\right\vert :t\in
A\right\} ,\sup \left\{ \left\vert g(\varphi (t))\right\vert :t\in B\right\}
\right)  \\
&=&\max \left( \left\Vert f\right\Vert ,\left\Vert g\right\Vert \right)
=\left\Vert (f,g)\right\Vert _{0}\text{.}
\end{eqnarray*}%
Therefore, it can be easily observed that $Q$ is an isometric isomorphism of 
$C(\omega ^{\ast })\times C(\omega ^{\ast })$ onto $C(\omega ^{\ast })$. $%
\square \medskip $

From the above theorem it can be proven that any nonreflexive Banach space
ultrapower $\ell _{\infty }(E)/\mathcal{N}_{\mathcal{U}}$ can be represented
as a direct sum of two closed subspaces which are isometrically isomorphic
to $\ell _{\infty }(E)/\mathcal{N}_{\mathcal{U}}$ $($cf. \cite{LW, M}$)$%
.\medskip 

\textbf{Corollary }$\mathbf{14}$ $(CH)$. \textit{Let }$E$\textit{\ be any
infinite dimensional nonsuperreflexive Banach space and} $\ell _{\infty }(E)/%
\mathcal{N}_{\mathcal{U}}$ \textit{be} \textit{its ultrapower. Then} $\ell
_{\infty }(E)/\mathcal{N}_{\mathcal{U}}=M\oplus N$ \textit{where} $M$ 
\textit{and} $N$ \textit{are closed subspaces of this ultrapower. Both
subspaces} $M$ \textit{and} $N$ \textit{are isometrically isomorphic to the
ultrapower} $\ell _{\infty }(E)/\mathcal{N}_{\mathcal{U}}$.

\textit{Proof}. Suppose that $Q$ is the operator from the proof of the
previous theorem mapping $C(\omega ^{\ast })\times C(\omega ^{\ast })$ onto $%
C(\omega ^{\ast })$. Then these subspaces $M$ and $N$ can be represented as $%
M=Q(\ell _{\infty }(E)/\mathcal{N}_{\mathcal{U}}\times \{0\})$ and $%
N=Q(\{0\}\times \ell _{\infty }(E)/\mathcal{N}_{\mathcal{U}})$.
Consequently, $M$ and $N$ are isometrically isomorphic to the spaces $\ell
_{\infty }(E)/\mathcal{N}_{\mathcal{U}}\times \{0\}$ and $\{0\}\times \ell
_{\infty }(E)/\mathcal{N}_{\mathcal{U}}$ $($respectively$)$ and each of
these spaces is isometrically isomorphic to the Banach space ultrapower $%
\ell _{\infty }(E)/\mathcal{N}_{\mathcal{U}}$. From the fact that $M$ and $N$
are closed subspaces it follows that $\ell _{\infty }(E)/\mathcal{N}_{%
\mathcal{U}}\times \{0\}$ and $\{0\}\times \ell _{\infty }(E)/\mathcal{N}_{%
\mathcal{U}}$ are complementary subspaces. Hence, it can be proved that $%
\ell _{\infty }(E)/\mathcal{N}_{\mathcal{U}}=Q(\ell _{\infty }(E)/\mathcal{N}%
_{\mathcal{U}}\times \ell _{\infty }(E)/\mathcal{N}_{\mathcal{U}})=Q(\ell
_{\infty }(E)/\mathcal{N}_{\mathcal{U}}\times \{0\})\oplus Q(\{0\}\times
\ell _{\infty }(E)/\mathcal{N}_{\mathcal{U}})=M\oplus N$. Therefore, the
following formula is obtained $M\cong N\cong \ell _{\infty }(E)/\mathcal{N}_{%
\mathcal{U}}$. $\square \bigskip $

$\mathbf{6.}$\textbf{\ Concluding remarks\bigskip }

It should be stressed that the relation of isometric isomorphism between
nonreflexive Banach space ultrapowers and the space of continuous, bounded
and real-valued functions on the Parovicenko space $($which occurs under $%
CH) $, i.e., $\ell _{\infty }(E)/\mathcal{N}_{\mathcal{U}}\cong C(\omega
^{\ast })$ is very helpful in proving theorems about these model-theoretical
constructions. It was demonstrated that if the \textit{Continuum Hypothesis}
holds and $E$ is any infinite dimensional nonsuperreflexive Banach space,
then its ultrapower $\ell _{\infty }(E)/\mathcal{N}_{\mathcal{U}}$ is always
congruent to the space $C(\omega ^{\ast })$. This identification allowed to
answer several questions concerning the structure of nonreflexive Banach
space ultrapowers.

In the future it is planned to show $($using the above mentioned
Representation Theorem$)$ that all nonreflexive Banach space ultrapowers are 
$($as Banach spaces$)$ \textit{primary}. Consequently, this fact will enable
to demonstrate that every nonsuperreflexive Banach space can be represented
- from the model-theoretical point of view - in the form of its primary
ultrapower.

\bigskip

Piotr WILCZEK, Foundational Studies Center, ul. Na Skarpie 99/24, 61-163
Pozna\'{n}, POLAND, edwil@mail.icpnet.pl \bigskip

\end{document}